\documentclass[12pt]{article}
\usepackage{latexsym,amsfonts,amssymb}
\usepackage[dvips]{color}
\usepackage{colordvi,multicol}

\def \cal{\mathcal}

\newtheorem{thm}{Theorem}[section]

\newtheorem{pro}[thm]{Proposition}

\newtheorem{rem}[thm]{Remark}

\begin{document}
\noindent {\small Stochastic Processes and their Applications, 122 (2012) 2319-2328}

\vskip 0.6cm

%\title{\bf Hunt's Hypothesis (H) and Getoor's Conjecture for L\'{e}vy Processes}
%\author{ Ze-Chun Hu and Wei Sun} \maketitle
%\date{}

\centerline{\Large\bf Hunt's hypothesis (H) and Getoor's conjecture}
\centerline{\Large\bf for L\'{e}vy processes}

\vskip 1cm \centerline{Ze-Chun Hu} \centerline{\small Department
of Mathematics, Nanjing University, Nanjing 210093, China}
 \centerline{\small e-mail: huzc@nju.edu.cn}
\vskip 0.7cm \centerline{Wei Sun} \centerline{\small Department of
Mathematics and Statistics, Concordia University,}
\centerline{\small Montreal, H3G 1M8, Canada} \centerline{\small
e-mail: wsun@mathstat.concordia.ca} \vskip 1cm

\noindent{\bf Abstract}\quad In this paper, Hunt's hypothesis (H)
and Getoor's conjecture for L\'{e}vy processes are revisited. Let
$X$ be a L\'{e}vy process  on $\mathbf{R}^n$ with
L\'{e}vy-Khintchine exponent $(a,A,\mu)$. {First, we show that if
$A$ is non-degenerate then $X$ satisfies (H). Second, under the
assumption that $\mu({\mathbf{R}^n\backslash
\sqrt{A}\mathbf{R}^n})<\infty$,  we show that $X$ satisfies (H) if
and only if  the equation
$$
\sqrt{A}y=-a-\int_{\{x\in {\mathbf{R}^n\backslash
\sqrt{A}\mathbf{R}^n}:\,|x|<1\}}x\mu(dx),\ y\in \mathbf{R}^n,
$$
has at least one solution. Finally, we show that if $X$ is a subordinator and satisfies (H) then its drift coefficient must be 0.}

\smallskip


\smallskip


%\tableofcontents

\section{Introduction and main results}
Let $X$ be a nice Markov process. Hunt's hypothesis (H) says that
``every semipolar set of $X$ is polar".  (H) plays a crucial role in the potential theory of
(dual) Markov processes. We refer the reader to Blumenthal and Getoor \cite[Chapter
VI]{BG68}, \cite{BG70}  for details. In spite of
its importance, (H) has been verified only in some special
situations.  Let
$X$ and $\hat{X}$ be a pair of dual Markov processes as in \cite[Chapter
VI]{BG68}. Then, (H) holds if and only if the fine and cofine topologies differ by polar sets, see \cite[VI.4.10]{BG68} and Glover \cite[Theorem (2.2)]{G83}.
Some forty years ago, Getoor conjectured that essentially  all
L\'{e}vy processes {satisfy (H).}

Throughout this paper, we let $(\Omega,{\cal F},P)$ be a probability space and $X=(X_t)_{t\ge 0}$ be a $\mathbf{R}^n$-valued L\'{e}vy process on $(\Omega,{\cal F},P)$ with L\'{e}vy-Khintchine exponent
$\psi$, i.e.,
\begin{eqnarray*}
E[\exp\{i\langle z,X_t\rangle\}]=\exp\{-t\psi(z)\},\  z\in
\mathbf{R}^n,t\ge 0,
\end{eqnarray*}
where $E$ denotes the expectation {with respect to} $P$. For $\psi$, we have the following famous
L\'{e}vy-Khintchine formula:
\begin{eqnarray*}
\psi(z)=i\langle a,z\rangle+\frac{1}{2}\langle z,Az\rangle+\int_{\mathbf{R}^n}
\left(1-e^{i\langle z,x\rangle}+i\langle z,x\rangle 1_{\{|x|<1\}}\right)\mu(dx),
\end{eqnarray*}
where $a\in \mathbf{R}^n,A$ is a symmetric nonnegative definite $n\times n$ matrix, and
$\mu$ is a measure (called the L\'evy measure) on $\mathbf{R}^n\backslash\{0\}$
satisfying $\int_{\mathbf{R}^n\backslash\{0\}} (1\wedge |x|^2)\mu(dx)<\infty$. Hereafter, we use Re$(\psi)$ and
Im$(\psi)$ to denote the real part and imaginary part of $\psi$, respectively, and
use $(a,A,\mu)$ to denote $\psi$ sometimes. For every $x\in \mathbf{R}^n$, we denote by $P^x$ the law of $x+X$ under $P$. In particular, $P^0=P$.

Let $B\subset \mathbf{R}^n$. We define the first hitting time of $B$ by
$$
\sigma_B:=\inf\{t>0:X_t\in B\}.
$$
Denote by ${\cal B}^*$ the family of all nearly Borel sets relative
to $X$ (cf. \cite[I.10.21]{BG68}). A set $B\subset \mathbf{R}^n$ is
called polar (resp. essentially polar) if there exists a set $C\in
{\cal B}^*$ such that $B\subset C$ and $P^x(\sigma_C<\infty)=0$ for
every $x\in \mathbf{R}^n$ (resp. $dx$-almost every $x\in
\mathbf{R}^n$). Hereafter $dx$ denotes the Lebesgue measure on
$\mathbf{R}^n$. $B$ is called a thin set if there exists a set $C\in
{\cal B}^*$ such that $B\subset C$ and {$P^x(\sigma_C=0)=0$} for
every $x\in \mathbf{R}^n$. $B$ is called semipolar if
$B\subset\cup_{n=1}^{\infty}B_n$ for some thin sets
$\{B_n\}_{n=1}^{\infty}$.

Before introducing our results, we first recall some important results obtained so far for Getoor's conjecture. When $n=1$, Kesten \cite{Ke69} {(cf. also Bretagnolle \cite{Br71}) showed   that if $X$ is not a compound Poisson process, then  every $\{x\}$ is non-polar} if and only if
\begin{eqnarray*}\label{Ke69-a}
\int_0^{\infty}\mbox{Re}([1+\psi(z)]^{-1})dz<\infty.
\end{eqnarray*}
(If $X$ is a compound Poisson {process}, then it is easy to see that every $x$ {is regular} for $\{x\}$, i.e., $P^x(\sigma_{\{x\}}>0)=0$.)
Port and Stone \cite{PS69} proved that for the asymmetric
Cauchy process on the line every $x$ is regular for $\{x\}$. Hence only the empty set is a semipolar set and therefore (H) holds in this case. Further, Blumenthal and Getoor
\cite{BG70} showed that all  stable processes with index
$\alpha\in (0,2)$ on the line satisfy (H).

 Kanda \cite{Ka76} and Forst \cite{F75} proved that (H)
holds if $X$ has bounded continuous transition densities (with respect to $dx$) and
the L\'{e}vy-Khintchine exponent $\psi$ satisfies $|\mbox{Im} (\psi)|\leq
M(1+\mbox{Re}(\psi))$ for some positive constant $M$.
Rao \cite{R77} gave a short proof of the Kanda-Forst theorem under the weaker condition that $X$ has resolvent densities.
 In particular, for $n>1$ all
stable processes of index $\alpha\neq 1$ satisfy (H). Kanda
\cite{Ka78} settled this problem for the case $\alpha=1$
assuming the linear term vanishes. Silverstein {\cite{Si77}
extended the Kanda-Forst condition to the non-symmetric Dirichlet
forms setting, and  Fitzsimmons \cite{Fi01} extended it to the
semi-Dirichlet forms setting. Glover and Rao \cite{GR86}
proved that $\alpha$-subordinates of general Hunt processes satisfy
(H).
%Rao  \cite{R87}(1987) derived some of the above-mentional
%results, using some general energy considerations developed in
%Pop-Stojanovic and Rao \cite{PR81}(1981). By using the results in
%\cite{R87},
{Rao \cite{R88} proved that if all 1-excessive functions of $X$ are
lower semicontinuous and $|{\rm Im}(\psi)|\leq (1+{\rm
Re}(\psi))f(1+{\rm Re}(\psi))$, where $f$ is an increasing function
on $[1,\infty)$ such that $\int_N^{\infty}(zf(z))^{-1}dz=\infty$ for
any $N\geq 1$, then $X$ satisfies (H).}

Now we introduce the main results of this paper. To state the first
result, we let $\bar{X}$ be an independent copy of $X$. Define the
symmetrization $\tilde{X}$ of $X$ by $\tilde{X}:=X-\bar{X}$.

\begin{thm}\label{thm1.1}
Suppose that $A$ is non-degenerate, i.e.,
$A$ is of full rank. Then:
\begin{itemize}
\item[(i)] $X$ satisfies  (H);
\item[(ii)] The Kanda-Forst condition $|{\rm Im} (\psi)|\leq M(1+{\rm Re}(\psi))$ holds for some
positive constant $M$;
\item[(iii)] $X$  and $\tilde{X}$ have the same polar sets.
\end{itemize}
\end{thm}

Denote $b:=-a$ and $\mu_1:=\mu|_{\mathbf{R}^n\backslash
\sqrt{A}\mathbf{R}^n}$. If
$\int_{|x|<1}|x|\mu_1(dx)<\infty$, we set $b':=b-\int_{|x|<1}x\mu_1(dx)$. To state the second result, we define the following solution condition:

\noindent $(S)$\hspace{1cm} The equation $\sqrt{A}y=b',\ y\in \mathbf{R}^n$,
has at least one solution.

\begin{thm}\label{thm1.2}
Suppose that $\mu({\mathbf{R}^n\backslash
\sqrt{A}\mathbf{R}^n})<\infty$. Then, the following three claims are
equivalent:
\begin{itemize}
\item[(i)] $X$ satisfies  (H);
\item[(ii)]  (S) holds;
\item[(iii)] The Kanda-Forst condition $|{\rm Im} (\psi)|\leq M(1+{\rm Re}(\psi))$ holds for some
positive constant $M$.
\end{itemize}
\end{thm}

\begin{rem} (i) Theorem \ref{thm1.1} tells us that if a L\'evy process on $\mathbf{R}^n$ is perturbed by an independent (small) $n$-dimensional Brownian motion, then the perturbed L\'evy process must satisfy (H).

(ii) By Theorem \ref{thm1.2} and Jacob \cite[Example
4.7.32]{Ja01}, one finds that if $X$ satisfies  (H) and $\mu({\mathbf{R}^n\backslash
\sqrt{A}\mathbf{R}^n})<\infty$, then $X$ must be associated
with a Dirichlet form on $L^2(\mathbf{R}^n;dx)$.
\end{rem}

\begin{pro}\label{pro1.3}
Suppose that $\mu({\mathbf{R}^n\backslash
\sqrt{A}\mathbf{R}^n})<\infty$. Then:
\begin{itemize}
\item[(i)] $X$ has transition densities implies that all the claims of Theorem \ref{thm1.2} are fulfilled.
\item[(ii)] If one of the claims in Theorem \ref{thm1.2} is fulfilled, then the following four claims are
equivalent:
\begin{itemize}
\item[(a)] Every essentially polar set of $X$ is polar;
\item[(b)] $X$ has resolvent densities;
\item[(c)] $X$ has transition densities;
\item[(d)] $A$ is of full rank.
\end{itemize}
\end{itemize}
\end{pro}

\begin{pro}\label{pro1.6}
Suppose that $X$ has bounded continuous transition densities,
 and $X$ and $\tilde{X}$ have the same polar sets.
Then $X$ satisfies  (H).
\end{pro}

Suppose that $X$ is a subordinator. Then $\psi$ can be expressed by
$$
\psi(z)=-idz+\int_{(0,\infty)}{\left(1-e^{izx}\right)}\mu(dx),\
z\in \mathbf{R},
$$
where $d\geq 0$ (called the drift coefficient) and $\mu$ satisfies
$\int_{(0,\infty)}(1\wedge x)\mu(dx)<\infty$.

\begin{pro}\label{pro1.7}
If $X$ is a subordinator and satisfies (H), then $d=0$.
\end{pro}

The rest of this paper is organized as follows. In Section 2, we
recall the L\'{e}vy-It\^{o} decomposition of L\'{e}vy processes and
discuss the orthogonal transformation of L\'{e}vy processes. In
Section 3, we present the proofs of our results.

\section{L\'{e}vy-It\^{o} decomposition and orthogonal transformation of L\'{e}vy processes}
\setcounter{equation}{0}
\subsection{L\'{e}vy-It\^{o} decomposition}

%Firstly, let's recall the famous L\'{e}vy-It\^{o} decompositions for
%L\'{e}vy processes.

\begin{thm} {\bf (L\'{e}vy-It\^{o})}
Let $X$ be a L\'{e}vy process  on $\mathbf{R}^n$ with  {exponent}
$(a,A,\mu)$. Then there exist a Brownian motion $B_A$ on $\mathbf{R}^n$ with
covariance matrix $A$ and an independent Poisson random measure $N$
on $\mathbf{R}^+\times (\mathbf{R}^n\backslash\{0\})$ such that, for
each $t\geq 0$,
\begin{eqnarray}\label{IL}
X_t=bt+B_A(t)+\int_{|x|\geq 1}xN(t,dx)+\int_{|x|<1}x\tilde{N}(t,dx),
\end{eqnarray}
where $b=-a, \tilde{N}(t,F)=N(t,F)-t\mu(F)$.
\end{thm}

Define
\begin{eqnarray*}
&&X^{(I)}_t:=bt+B_A(t),\ X^{(II)}_t:=\int_{|x|\geq 1}xN(t,dx),\
X^{(III)}_t:=\int_{|x|<1}x\tilde{N}(t,dx),\ \ t\ge0.
\end{eqnarray*}
Then $X^{(I)}$, $X^{(II)}$ and $X^{(III)}$ are mutually
independent, $X^{(II)}$ is a compound Poisson process, and
$X^{(III)}$ is a square integrable martingale. For convenience,
we write $X^{(I)}_t=bt+\sqrt{A}B_t$, where $B=(B_t)_{t\ge 0}$ is a standard
$n$-dimensional Brownian motion.

\subsection{Orthogonal transformation}

Let $X$ be a L\'{e}vy process on $\mathbf{R}^n$ with
 exponent $(a,A,\mu)$. Since $A$ is a symmetric nonnegative definite matrix, there exists an orthogonal matrix
$O$ such that
$$
OAO^T=diag(\lambda_1,\dots,\lambda_n):=D,
$$
where $\lambda_1\geq \lambda_2\geq \cdots\geq \lambda_n\geq 0$ and $O^T$ denotes the transpose of $O$.
We fix such an orthogonal matrix
$O$ and define $ Y_t:=OX_t$, $t\geq 0$. Then $Y=(Y_t)_{t\ge 0}$ is a L\'{e}vy
process on $\mathbf{R}^n$  and $X$ satisfies (H) if and only if
$Y$ satisfies (H). We will see that sometimes it is more convenient to work with $Y$.
%First, we derive the  exponent of $Y$. By the expression of the
% exponent of $X$, we get
%\begin{eqnarray*}
%&&E[e^{i\langle z,Y_t\rangle}]=E[e^{i\langle z,OX_t\rangle}]=E[e^{i\langle O^Tz,X_t\rangle}]\\
%&&=\exp\left[-t\left(i\langle a,O^Tz\rangle+\frac{1}{2}\langle O^Tz,AO^Tz\rangle+\int_{\mathbf{R}^n}\left(1-e^{i\langle O^Tz,x\rangle}+i\langle O^Tz,x\rangle 1_{\{|x|<1\}}\right)\mu(dx)\right)\right]\\
%&&=\exp\left[-t\left(i\langle Oa,z\rangle+\frac{1}{2}\langle
%z,Dz\rangle+\int_{\mathbf{R}^n}\left(1-e^{i\langle z,Ox\rangle}+
%i\langle z,Ox\rangle 1_{\{|x|<1\}}\right)\mu(dx)\right)\right]\\
%&&=\exp\left[-t\left(i\langle Oa,z\rangle+\frac{1}{2}\langle
%z,Dz\rangle+\int_{\mathbf{R}^n}\left(1-e^{i\langle z,y\rangle}+
%i\langle z,y\rangle 1_{\{|y|<1\}}\right)\mu
%O^{-1}(dy)\right)\right],
%\end{eqnarray*}
By the expression of the  exponent of $X$ and simple computation, we get that the  exponent of $Y$ is
$(Oa,D,\mu O^{-1})$,
where $\mu O^{-1}(B)=\mu(\{x\in \mathbf{R}^n:Ox\in B\})$ for any
Borel set $B$ of $\mathbf{R}^n$.

From now on, we denote by $k$ the rank of $A$. Then, the orthogonal transformation satisfies the following properties:

\noindent (1)  $\mu({\mathbf{R}^n\backslash
\sqrt{A}\mathbf{R}^n})<\infty$ if and only if  $\mu O^{-1}$ is  a
finite measure on $\mathbf{R}^k\times (\mathbf{R}^{n-k}\backslash
\{0\})$. (When $k=n$, ${\mathbf{R}^n\backslash
\sqrt{A}\mathbf{R}^n}$ and $\mathbf{R}^k\times (\mathbf{R}^{n-k}\backslash
\{0\})$ are the empty set.)

\noindent (2) If $\int_{|x|<1}|x|\mu_1(dx)<\infty$, then
$$
\int_{\{y\in \mathbf{R}^k\times
(\mathbf{R}^{n-k}\backslash\{0\}):|y|<1\}}|y|\mu O^{-1}(dy)=\int_{|y|<1}|y|\mu_1 O^{-1}(dy)=\int_{|x|<1}|x|\mu_1(dx)<\infty.
$$
Recall that $b'=b-\int_{|x|<1}x\mu_1(dx)$. Define $\bar{b}:=Ob'$. Then
$$
\bar{b}=Ob-\int_{\{y\in \mathbf{R}^k\times
(\mathbf{R}^{n-k}\backslash\{0\}):|y|<1\}}y\,\mu O^{-1}(dy).
$$
Note that $\sqrt{A}=O^T\sqrt{D}O$. Then, the equation $\sqrt{A}y=b'$ is equivalent to $\sqrt{D}Oy=Ob'$.
Therefore,  the equation $\sqrt{A}y=b', y\in \mathbf{R}^n$, has a solution if and
only if the equation $\sqrt{D}y=\bar{b}, y\in \mathbf{R}^n$, has a solution.

 \noindent (3) Suppose that $\int_{|x|<1}|x|\mu_1(dx)<\infty$. Then, {by the L\'{e}vy-It\^{o} decomposition (\ref{IL}),
$Y$ can be expressed by
\begin{eqnarray}\label{3.2-a}
Y_t=Obt+\sqrt{D}\bar{B}_t+\int_{|y|\geq
1}y\bar{N}(t,dy)+\int_{|y|<1}y\tilde{\bar{N}}(t,dy),
\end{eqnarray}
where $\bar{B}=OB$ is a standard Brownian motion on $\mathbf{R}^n$, $\bar{N}$ is a Poisson random measure on $\mathbf{R}^+\times (\mathbf{R}^n\backslash \{0\})$ with $\mu O^{-1}$ being its intensity measure, $
\tilde{\bar{N}}(t,F)=\bar{N}(t,F)-t\,\mu O^{-1}(F)$, $\bar{B}$ and $\bar{N}$ are independent. We rewrite
(\ref{3.2-a}) as
\begin{eqnarray*}
Y_t=Y^{(1)}_t+Y^{(2)}_t,
\end{eqnarray*}
where
\begin{eqnarray*}
&&Y^{(1)}_t:=\bar{b}t+\sqrt{D}\bar{B}_t+\int_{\{y\in
\mathbf{R}^k\times\{0\}:|y|\geq 1\}}y\bar{N}(t,dy)+\int_{\{y\in
\mathbf{R}^k\times\{0\}:|y|< 1\}}y\tilde{\bar{N}}(t,dy),\\
&&Y^{(2)}_t:=\int_{\mathbf{R}^k\times (\mathbf{R}^{n-k}\backslash
\{0\})}y\bar{N}(t,dy),\nonumber
\end{eqnarray*}
 $Y^{(1)}$ and $Y^{(2)}$ are independent.

By (2), we can see that (S) {holds if and only if} $\bar{b}\in
\mathbf{R}^k\times \{0\}$. In this case, $Y^{(1)}$ can be regarded
as a $k$-dimensional L\'{e}vy process on $\mathbf{R}^k\times \{0\}$, which has a non-degenerate
Gaussian component. If $\mu(\mathbf{R}^n\backslash
\sqrt{A}\mathbf{R}^n)<\infty$, then  $Y^{(2)}$ is a compound Poisson
process.}

\section{Proofs of the main results}\setcounter{equation}{0}
 \subsection{Proof of Theorem \ref{thm1.1}}

First, we  prove (ii). Since $A$ is of full rank,
there exists a constant $c>0$ such that $\langle z,Az\rangle\ge
c\langle z,z\rangle$, $\forall z\in \mathbf{R}^n$. Then
\begin{eqnarray}\label{proof-thm1.1-a}
{\rm Re}\psi(z)=\frac{1}{2}\langle
z,Az\rangle+\int_{\mathbf{R}^n}(1-\cos\langle z,x\rangle)\mu(dx)\geq
\frac{1}{2}\langle z,Az\rangle\geq \frac{c}{2}\langle z,z\rangle.
\end{eqnarray}
By the Cauchy-Schwarz inequality, one finds that
 $|\langle a,z\rangle|$ is
controlled by $1+{\rm Re}\psi(z)$. To establish the Kanda-Forst
condition, we need only show that $|{\rm
Im}\{\int_{|x|<1}\left(1-e^{i\langle z,x\rangle}+i\langle
z,x\rangle\right)\mu(dx)\}|$ is controlled by $\langle z,z\rangle$.
Note that $|t-\sin t|\le t^2/2$ for any $t\in \mathbf{R}$. Then,
\begin{eqnarray*}
\left|{\rm Im}\left\{\int_{|x|<1}\left(1-e^{i\langle z,x\rangle}+i\langle z,x\rangle\right)\mu(dx)\right\}\right|
&=&\left|\int_{|x|<1}\left(\langle z,x\rangle-\sin\langle z,x\rangle\right)\mu(dx)\right|\\
&\le&\frac{1}{2}\int_{|x|<1}|\langle z,x\rangle|^2\mu(dx)\\
&\le&\left(\frac{1}{2}\int_{|x|<1}|x|^2\mu(dx)\right)|z|^2.
\end{eqnarray*}
Therefore (ii) holds.

Second, we prove (i). By (\ref{proof-thm1.1-a}), we get
\begin{eqnarray}\label{proof-thm1.1-b}
\lim_{|z|\to\infty}\frac{{\rm Re}\psi(z)}{\ln(1+|z|)}=\infty.
\end{eqnarray}
By Hartman and Wintner \cite{HW42} (cf. also Knopova  and Schilling
\cite{KS10}) and (\ref{proof-thm1.1-b}), we find that $X$ has
bounded continuous transition densities. Then, by (ii) and  the Kanda-Forst theorem,
we obtain (i).

Finally, we prove (iii). Denote by $\tilde{\psi}$ the
L\'{e}vy-Khintchine exponent of $\tilde{X}$. Note that
$\tilde{\psi}=2{\rm Re}(\psi)$. Then, for any $\lambda\geq 1$, by
(ii) we get (cf. Kanda \cite[Page 163]{Ka76})
\begin{eqnarray}\label{proof-thm1.1-c}
2{\rm Re}\left(\frac{1}{\lambda+\tilde{\psi}(\xi)}\right)
&=&\frac{1}{\frac{1}{2}\lambda+\frac{1}{2}\tilde{\psi}(\xi)}\geq
\frac{1}{\lambda+{\rm Re}\psi(\xi)}\nonumber\\
 &\geq&{\rm
Re}\left(\frac{1}{\lambda+\psi(\xi)}\right)\nonumber\\
&=&\frac{\lambda+{\rm Re}\psi(\xi)}{(\lambda+{\rm Re}\psi(\xi))^2+({\rm Im}\psi(\xi))^2}\nonumber\\
&=&\frac{1}{\lambda+{\rm Re}\psi(\xi)}\left[1+\left(\frac{{\rm
Im}\psi(\xi)}{\lambda+{\rm
Re}\psi(\xi)}\right)^2\right]^{-1}\nonumber\\
&\geq&\frac{1}{\lambda+{\rm
Re}\psi(\xi)}\left[1+\left(\frac{M(1+{\rm
Re}\psi(\xi))}{\lambda+{\rm
Re}\psi(\xi)}\right)^2\right]^{-1}\nonumber\\
&\geq&\frac{1}{(1+M^2)(\lambda+\frac{1}{2}\tilde{\psi}(\xi))}\nonumber\\
&\geq&\frac{1}{1+M^2}{\rm
Re}\left(\frac{1}{\lambda+\tilde{\psi}(\xi)}\right).
\end{eqnarray}
By (\ref{proof-thm1.1-c}), { the above proved fact that $X$
has bounded continuous transition densities, and Kanda \cite[Theorem 1]{Ka76} (or Hawkes \cite[Theorems 2.1
and 3.3]{Ha79}), we obtain (iii).}\hfill\fbox

\subsection{Proof of Theorem \ref{thm1.2}}
By the discussion of \S2.2, we know that $X$ satisfies (H) if and only if
$Y$ satisfies (H), and (S) holds if and only if $\bar{b}\in
\mathbf{R}^k\times \{0\}$. By the expression of the  exponent of $Y$, it is easy to see that the Kanda-Forst condition holds for $X$ if and only if it holds for $Y$. Hence, to prove Theorem \ref{thm1.2}, we may and do assume without loss of generality that {$A=diag(\lambda_1,\ldots,\lambda_n):=D$, where $\lambda_1\geq \lambda_2\geq \cdots\geq \lambda_k>0,\lambda_{k+1}=\lambda_{k+2}=\cdots=\lambda_n=0\ (k\geq 0)$, and $X$ has the expression
\begin{eqnarray*}\label{3.2-b}
X_t=X^{(1)}_t+X^{(2)}_t,\ \ t\ge0,
\end{eqnarray*}
where
\begin{eqnarray}\label{test}
&&X^{(1)}_t:=b't+\sqrt{D}B_t+\int_{\{x\in
\mathbf{R}^k\times\{0\}:|x|\geq 1\}}xN(t,dx)+\int_{\{x\in
\mathbf{R}^k\times\{0\}:|x|< 1\}}x\tilde{N}(t,dx),\\
&&X^{(2)}_t:=\int_{\mathbf{R}^k\times (\mathbf{R}^{n-k}\backslash
\{0\})}xN(t,dx),\nonumber
\end{eqnarray}
$b'$ is the same as in \S1, and $B$, $N$ and $\tilde{N}$ are the same as in \S2.1. }

If $k=0$, then $X_t=b't+X^{(2)}_t$. Since $X^{(2)}$ is a compound Poisson process, it
is easy to see that (i), (ii) and (iii) are equivalent in this case. Below we assume that $k\ge 1$.

\noindent (ii) $\Rightarrow$ (iii): Suppose that (S) holds, i.e.,
$b'\in \mathbf{R}^k\times \{0\}$. Then $X^{(1)}$ stays in
$\mathbf{R}^k\times \{0\}$ if it starts there. By Theorem
\ref{thm1.1}, the Kanda-Forst condition
 holds for $X^{(1)}$. Since $X^{(2)}$ is a compound Poisson process,
 its L\'evy-Khintchine exponent is bounded. Hence the Kanda-Forst
 condition holds for $X$, i.e., (iii) holds.

\bigskip

\noindent (iii) $\Rightarrow$ (ii): Suppose that the Kanda-Forst
 condition holds for $X$. Since the L\'evy-Khintchine exponent of
 $X^{(2)}$ is  bounded, we get that the Kanda-Forst condition
 holds for $X^{(1)}$. Assume that $b'\notin
\mathbf{R}^k\times \{0\}$. We will reach a contradiction. Denote
$b'=(b'_1,\ldots,b'_n)$. Without loss of generality, we assume that
$b'_n\neq 0$.
 Let $\psi_1$ be the L\'evy-Khintchine exponent of $X^{(1)}$. Then
\begin{eqnarray*}
\psi_1(z)=i\langle b',z\rangle+\frac{1}{2}\langle
z,\sqrt{D}z\rangle+\int_{\mathbf{R}^k\times \{0\}}
\left(1-e^{i\langle z,x\rangle}+i\langle z,x\rangle
1_{\{|x|<1\}}\right)\mu(dx).
\end{eqnarray*}
It follows that if $z=(z_1,\ldots,z_n)$ with $z_i=0$, $i=1,\ldots,n-1$
and $z_n\neq 0$, then $\psi_1(z)=b_n'z_ni$ and thus the Kanda-Forst condition cannot hold for $X^{(1)}$. Hence  $b'\in
\mathbf{R}^k\times \{0\}$ and therefore (S) holds.

\bigskip

\noindent (i) $\Rightarrow$ (ii):  We will show $b'\notin \mathbf{R}^k\times \{0\}$ implies that $X$ does not satisfy (H). We first consider the case that $\mu_1\not= 0$.

Suppose that $b'\notin \mathbf{R}^k\times \{0\}$. First, we
show that $\mathbf{R}^k\times \{0\}$ is a thin set of $X$.   Let
$T^{(2)}_1$ be
 the first jumping time of $X^{(2)}$. Since $X^{(2)}$ is a compound Poisson process,  $T^{(2)}_1$
 has an exponential distribution, in particular,
\begin{eqnarray}\label{thm1.1-b}
P(T^{(2)}_1>0)=1.
\end{eqnarray}
For any $x\in
\mathbf{R}^k\times \{0\}$ and any $t>0$, we know that
$x+X^{(1)}_t\notin \mathbf{R}^k\times \{0\}$ since $b'\notin
\mathbf{R}^k\times \{0\}$, which together with (\ref{thm1.1-b})
implies that
 \begin{eqnarray}\label{thm1.1-c}
{P^{x}(\sigma_{\mathbf{R}^k\times \{0\}}=0)}&\leq &P^0\left(\exists t\in (0,T^{(2)}_1)\ s.t.\ x+X_t\in \mathbf{R}^k\times \{0\}\right)\nonumber\\
&=&P^0\left(\exists t\in (0,T^{(2)}_1)\ s.t.\ x+X^{(1)}_t\in
\mathbf{R}^k\times \{0\}\right)\nonumber\\
&=&0.
\end{eqnarray}
For any $x\notin \mathbf{R}^k\times \{0\}$, the distance between $x$
and {the subspace} $\mathbf{R}^k\times \{0\}$ is strictly positive. By
(\ref{thm1.1-b}) and the right continuity of the sample path of
$X^{(1)}$, we get
\begin{eqnarray}\label{thm1.1-d}
{P^{x}(\sigma_{\mathbf{R}^k\times \{0\}}=0)}&=&P^0\left(\exists \{t_n,n\geq 1\}\subset (0,T^{(2)}_1)\ s.t.\ x+X_{t_n}\in \mathbf{R}^k\times \{0\},t_n\downarrow 0\right)\nonumber\\
&=&P^0\left(\exists \{t_n,n\geq 1\}\subset (0,T^{(2)}_1)\ s.t.\
x+X^{(1)}_{t_n}\in \mathbf{R}^k\times \{0\},t_n\downarrow 0\right)\nonumber\\
&=&0.
\end{eqnarray}
It follows from (\ref{thm1.1-c}) and (\ref{thm1.1-d}) that $\mathbf{R}^k\times \{0\}$ is a
thin set  and thus a semipolar set of $X$.

{Next, we show that $\mathbf{R}^k\times \{0\}$ is not a polar set of $X$. Note that $P^0(T^{(2)}_1>s)>0$ for any $s>0$. Then{
\begin{eqnarray*}
P^{-b's}(\sigma_{\mathbf{R}^k\times \{0\}}<\infty)
&=&P^{-b's}\left(\exists\ t> 0\ s.t.\ X_t\in \mathbf{R}^k\times \{0\}\right)\\
&\ge&P^{-b's}\left(X_s\in \mathbf{R}^k\times \{0\}\right)\\
&=&P^0\left((X_s-b's)\in \mathbf{R}^k\times \{0\}\right)\\
&\ge&P^{0}(T^{(2)}_1>s)\\
&>&0.
\end{eqnarray*}}
Hence $\mathbf{R}^k\times \{0\}$ is not a polar set of $X$. Therefore
$X$ does not satisfy (H).

The case that $\mu_1=0$ can be proved  similarly by $T^{(2)}_1\equiv\infty$.

\bigskip

\noindent (ii) $\Rightarrow$ (i): Suppose that (S) holds, i.e., $b'\in
\mathbf{R}^k\times \{0\}$. Let $F$ be a semipolar set of $X$. We will show that $F$ is a
polar set of $X$. Without loss of generality, we assume that $F$ is a nearly Borel set. For $y\in
\mathbf{R}^{n-k}$, we define
$$
F_y:=F\cap (\mathbf{R}^k\times \{y\}).
$$
Since $X^{(2)}$ is a compound Poisson
process, one finds that $F_y$ is semipolar for the $k$-dimensional L\'evy process $(X^{(1)}, P^{(x,y)})_{x\in \mathbf{R}^{k}}$ on $\mathbf{R}^k\times \{y\}$. Hence $F_y$ is polar for
$(X^{(1)}, P^{(x,y)})_{x\in \mathbf{R}^{k}}$ by Theorem \ref{thm1.1}. Therefore,
\begin{equation}\label{apend1}
P\left(\exists t>0\ s.t.\ (x,y)+X^{(1)}_t\in F_y\right)=0,\ \ \forall x\in \mathbf{R}^{k}, \forall y\in \mathbf{R}^{n-k}.
\end{equation}

Denote  by
$\eta$ the distribution of $T^{(2)}_1$ under $P$. Let $\xi$ be a random variable taking values on $\mathbf{R}^k\times (\mathbf{R}^{n-k}\backslash
\{0\})$, which has distribution $\mu_1$ and is independent of
$X^{(1)}$ and $T^{(2)}_1$. Then, {for any} $x_0=(u,v)\in \mathbf{R}^{k}\times \mathbf{R}^{n-k}$, we obtain by (\ref{apend1}) that
\begin{eqnarray*}
P\left(x_0+X^{(1)}_{T^{(2)}_1}+\xi\in F\right)
&=& \int_{\mathbf{R}^k\times (\mathbf{R}^{n-k}\backslash
\{0\})}\int_{(0,\infty)}P((u,v)+X^{(1)}_t+(x,y)\in F)\eta(dt)
\mu_1(dx,dy)\\
&=& \int_{\mathbf{R}^k\times (\mathbf{R}^{n-k}\backslash
\{0\})}\int_{(0,\infty)}P((u+x,v+y)+X^{(1)}_t\in F_{v+y})\eta(dt)
\mu_1(dx,dy)\\
&=&0.
\end{eqnarray*}
Since $x_0$ is arbitrary, by the strong Markov
property of L\'evy process, $F$ is a polar set of $X$. Therefore, $X$ satisfies  (H).\hfill\fbox

\subsection{Proof of Proposition \ref{pro1.3}}

\noindent (i) Suppose that $X$ has transition densities. We will show that $A$ is of full rank. We adopt the setting of \S3.2. Assume that $k<n$. Set $X=(X^1,\dots,X^n)$ and $b'=(b'^1,\dots,b'^n)$. Without loss of generality, we suppose $\mu_1\neq
0$.  Let $T^{(2)}_1$ be
 the first jumping time of $X^{(2)}$.
 Then $T^{(2)}_1$ has an exponential distribution and thus
$P(T^{(2)}_1>1)>0$. It follows from (\ref{test}) that
$P(X_1^n=b'^n)>0$.
This contradicts with the assumption that $X$ has transition densities. Hence $A$ is of full rank. Therefore, the proof is completed by Theorem \ref{thm1.1}.

\bigskip

\noindent (ii) (a) $\Leftrightarrow$ (b) follows from
Fukushima \cite[(viii)]{Fu74}. (d) $\Rightarrow$ (c) $\Rightarrow$ (b) is
easy. (b) $\Rightarrow$ (c) follows from the Kanda-Forst condition,
Silverstein \cite[Theorem 3.2]{Si78}  and the spatial homogeneity of
L\'{e}vy processes. (Note that if the  L\'{e}vy process $X$ is associated with a Dirichlet form on $L^2(\mathbf{R}^n;dx)$, then the Dirichlet form is regular by Silverstein \cite[Lemma 1.5]{Si77}.) (c) $\Rightarrow$ (d) follows from the above
proof of (i). \hfill\fbox

\subsection{Proof of Proposition \ref{pro1.6}}
{The main idea has been used in the proof of Kanda
\cite[Theorem 2]{Ka78}. } Denote by $\tilde{\psi}$ the L\'{e}vy-Khintchine exponent of
$\tilde{X}$. Then, for any
$\lambda>0$,  we have
\begin{eqnarray*}\label{thm2.1-a}
{\rm Re}\left(\frac{1}{\lambda+\psi(\xi)}\right)\leq
\frac{1}{\lambda+{\rm Re}\psi(\xi)}=
\frac{1}{\lambda+\frac{1}{2}\tilde{\psi}(\xi)}\leq 2\ {\rm
Re}\left(\frac{1}{\lambda+\tilde{\psi}(\xi)}\right).
\end{eqnarray*}
 By Kanda \cite[Remark 2.1]{Ka76} or Hawkes \cite[Theorem 3.3]{Ha79}, we find that
 there exists a positive constant $M$ such that for every $\lambda>0$ and every compact $K$,
\begin{eqnarray}\label{thm2.1-b}
C^{\lambda}(K)\geq M\tilde{C}^{\lambda}(K),
\end{eqnarray}
where $C^{\lambda}(K)$ (resp. $\tilde{C}^{\lambda}(K)$) is
$\lambda$-capacity of $K$ relative to $X$ (resp. $\tilde{X}$). Since
$\tilde{X}$ is a symmetric L\'{e}vy process with
bounded continuous transition densities, it satisfies
 (H),  i.e., every semipolar set of $\tilde{X}$ is a polar set of
$\tilde{X}$. By Kanda \cite[Theorem 1]{Ka78}, we get
\begin{eqnarray}\label{thm2.1-c}
\lim_{\lambda\uparrow\infty}\tilde{C}^{\lambda}(K)=\infty
\end{eqnarray}
for every non-polar compact set $K$ of $\tilde{X}$. (We remark that, more generally,
(H) implies   (\ref{thm2.1-c}) under the weaker condition that $\tilde{X}$ has resolvent densities, see Getoor \cite[Theorem (11.21)]{Ge90}.)
 By the
assumption, we find that every non-polar compact set $K$ of $X$ is a
non-polar compact set of $\tilde{X}$. Thus, by (\ref{thm2.1-b}) and
(\ref{thm2.1-c}), we get
\begin{eqnarray*}
\lim_{\lambda\uparrow\infty}C^{\lambda}(K)=\infty
\end{eqnarray*}
for every non-polar compact set $K$ of $X$. Then, by Kanda
\cite[Theorem 1]{Ka78} again, we obtain that every semipolar set of
$X$ is a polar set of $X$.\hfill\fbox

\subsection{Proof of Proposition \ref{pro1.7}}

 Suppose that $d>0$. Then $X$ is strictly increasing, which together with the right continuity of sample paths implies that singletons are thin and thus semipolar. By Kesten \cite{Ke69} or  Bretagnolle \cite{Br71}, we know that $X$ hits points with positive probability. Hence (H) cannot hold. Therefore we must have $d=0$. \hfill\fbox

\smallskip

{ \noindent {\bf\large Acknowledgments} \vskip 0.1cm  \noindent
We thank an Associate Editor and a referee for their valuable suggestions, which improved the presentation of this paper.
We are
grateful to the support of NNSFC (Grant No. 10801072) and NSERC (Grant No. 311945-2008).}

\end{document}